\documentclass[11pt,reqno,fleqn,a4paper]{article}
\usepackage{amsmath,amssymb,amsfonts,amsthm,enumitem,cite}
\usepackage{setspace}
\usepackage{graphicx}
\usepackage{multicol}
\usepackage{lipsum,lineno}
\usepackage{ragged2e}
\usepackage{multirow}
\usepackage{bm}
\usepackage[margin=1.0in]{geometry}
\usepackage[utf8]{inputenc}
\usepackage[T1]{fontenc}
\usepackage{mathtools}
\usepackage[font=small,labelfont=bf]{caption}
\usepackage[subrefformat=parens]{subcaption}
\usepackage{mwe}
\usepackage{mfirstuc}

\providecommand{\keywords}[1]{\textbf{\textbf{\\Keywords. }}#1}
\providecommand{\class}[1]{\textbf{\textbf{Mathematics Subject Classification (2010). }}#1}
\usepackage{cite}
\MFUnocap{are}
\MFUnocap{or}
\MFUnocap{etc}

\def\capmystringaux#1#2\relax{\uppercase{#1}\lowercase{#2}}
\usepackage{lineno,hyperref}

\title{}
\author{}
\DeclareSymbolFont{slenderlargesymbols}{OMX}{ccex}{m}{n}
\usepackage{amsmath,amssymb,amsfonts,calrsfs,bbm}
\usepackage{setspace}
\usepackage{graphicx}
\usepackage{amsthm}
\usepackage{multicol}
\usepackage{hyperref}
\usepackage{lipsum,lineno}
\usepackage{amsthm}
\usepackage{ragged2e}
\usepackage{multirow}
\usepackage{bm}
\usepackage{mathtools}
\usepackage{mwe}
\usepackage{breqn}
\usepackage{mfirstuc}
\MFUnocap{are}
\MFUnocap{or}
\MFUnocap{etc}

\DeclarePairedDelimiterX{\Set}[2]\{\}{%
	\, #1 \;\delimsize\vert\; #2 \,
}

\usepackage{amsfonts}
\usepackage{mathrsfs}
\newcommand{\diam}{\operatorname{diam}}
\newcommand{\lipp}{\operatorname{Lip}_p}

\newtheorem{theorem}{Theorem}[section]

\theoremstyle{definition}
\newtheorem{definition}[theorem]{Definition}
\theoremstyle{remark}
\newtheorem{remark}[theorem]{Remark}
\newtheorem*{example}{Example}
\numberwithin{equation}{section}
\date{}
\begin{document}
	\title{\textbf{On iterated function systems and algebraic properties of Lipschitz maps in partial metric spaces}}
	\author{Praveen M\footnote{E-mail:praveenmaniyatt@gmail.com}, Sunil Mathew\footnote{E-mail:sm@nitc.ac.in}\\ \small Department of Mathematics,\\\small National Institute of Technology Calicut, Calicut - 673 601\\\small India.
}
	\maketitle	
\begin{abstract}
	This paper discusses, certain algebraic, analytic, and topological results on partial iterated function systems($IFS_p$'s). Also, the article proves the Collage theorem for partial iterated function systems. Further, it provides a method to address the points in the attractor of a partial iterated function system and obtain results related to the address of points in the attractor. The completeness of the partial metric space of contractions with a fixed contractivity factor is proved, under suitable conditions. Also, it demonstrates the continuity of the map that associates each contraction in a complete partial metric space to its corresponding unique fixed point. Further, it defines the $IFS_p$ semigroup and shows that under function composition, the set of Lipschitz transformations and the set of contractions are semigroups.\\
%
	\noindent\class{Primary 28A80; Secondary
		11B05.}\medskip
\end{abstract}

\section{Introduction}

\textit{Fractal theory} is a branch of mathematics emerged in the late $1960$'s which provides a novel geometrical view of the objects and phenomena in nature like coastlines, blood vessels, and pulmonary vessels, DNA, ocean waves, mountain ranges, heart sounds, river networks, trees, Romanesco broccoli, soil pores, snowflakes, etc. In recent years the fractal theory has achieved significant developments in both theoretical and applicable mathematics. The geometry proposed by the theory of fractals rectifies many flaws of the classical Euclidean geometry while modeling the natural objects and phenomena. Most of the early contributions to the fractal theory are credited to Benoit B. Mandelbrot, who coined the word \textit{fractal} in $1975$. 
The word fractal is derived from the Latin word \textit{fr\={a}ctus} for "broken" or "fractured". 
  Mandelbrot proposed a theory of \textit{roughness and self-similarity} in nature and demonstrated how we could obtain visual complexity from some simple rules. In his perspective, the things we consider as \textit{rough}, \textit{complex} or \textit{chaotic}, have an inherent \textit{degree of order}.

There are different methods available to create fractals such as iterated function systems or $IFSs$, strange attractors, L-systems, escape-time fractals, random fractals, finite subdivision rules, etc. The $IFS$ theory was introduced in $1981$ by John E Hutchinson\cite{hutchinson1979fractals}. The fractals constructed with $IFSs$ are mostly \textit{self-similar}, i.e., they are similar to a part of itself. An $IFS$ is defined formally as a finite set of contraction mappings on a complete metric space. There have been many generalizations to this theory over the past few years. The generalization in  which an infinite set of contractions is used in place of a finite set is called an \textit{infinite iterated function system}($IIFS$)\cite{fernau1994infinite} and is called a \textit{countable iterated function system}($CIFS$)\cite{secelean2013countable,secelean2012existence} when the set of contractions is countable. Also, there are generalizations in the literature with the contraction condition relaxed to generalized contraction conditions such as $r$-contraction\cite{secelean2013iterated}, convex contraction\cite{istratescu1981some,istractescu1982some,miculescu2015generalization,georgescu2017ifss}, Meir-Keeler type contraction, $F$-contraction\cite{secelean2013iterated}, weak contraction\cite{hata1985structure}, etc. The other types of generalization include the relaxation on the completeness condition of the underlying metric space, replacing the metric space with a product of metric spaces, relaxing the metric condition in the space to partial metric\cite{minirani2014fractals} condition, 
replacing the metric space with a general topological space\cite{mihail2012topological}, etc.

The experiences from the field of computer science motivate the notion of the partial metric, and its theory is an extension of the theory of metric spaces, with the significant difference being the self-distance can be non-zero in a partial metric space.
The theory of partial metric spaces was introduced in $1992$ by Steve G Matthews in connection with a study on dataflow networks. He also extended the Banach contraction principle to the partial metric spaces. In Matthews' definition of partial metric the range is restricted to $\mathbb{R}^+$ and S. J. O'Neill\cite{o1996partial} in $1996$ extended this range to $\mathbb{R}$. Even though  Matthews considered non-zero self distances, he restricted them to be small with the use of an additional condition in the definition of the partial metric. In $1997$, Reinhold Heckmann\cite{heckmann1999approximation} relaxed this condition to define the weak partial metric, in which large self-distances are possible. 

In this paper, we are providing certain algebraical, analytical and topological results concerning the partial iterated function systems($IFS_ps$), which was introduced in $2014$ by S. Minirani and Sunil Mathew\cite{minirani2014fractals} as a generalization of the $IFSs$ by changing the metric space in the system to a partial metric space.


\section{Preliminaries}
This section provides some of the basic definitions and results for the development of this paper which are from \cite{barnsley2014fractals,barnsley2006superfractals,matthews1994partial,minirani2014fractals,aydi2012partial,oltra2004banach,valero2005banach,bukatin2009partial,altun2010generalized}. Throughout this paper, $(X, d)$ denotes a complete metric space, $(X, p)$ denotes a complete partial metric space and $\mathbbm{H}(X)$ and $\mathbbm{H}_p(X)$ denotes the non-empty compact subsets of $(X, d)$ and $(X, p)$ respectively. 

The following definitions from literature helps in defining a metric called the Hausdorff metric in $\mathbbm{H}(X)$.
\begin{definition}\cite{barnsley2014fractals}
	Let $(X,d)$ be a metric space, $x\in X$ and $K,A,B\in \mathbbm{H}(X)$. Then $$d(x,K)=\min\{d(x,y) : y\in K\} \text{ and}$$
	$$d(A,B)=\max\{d(x,B) : x\in A\}.$$
\end{definition}
The reason that $d:\mathbbm{H}(X)\times \mathbbm{H}(X) \to \mathbb{R} $ is not a metric is explained in the following remark.
\begin{remark}\cite{barnsley2014fractals}
	In general, $d(A,B)\ne d(B,A)$ and $d(A,B)=0$ even if $A\ne B$.
\end{remark}
The Hausdorff metric on $\mathbbm{H}(X)$ is defined as follows.
\begin{definition}\cite{barnsley2014fractals} The map $h:\mathbbm{H}(X)\times \mathbbm{H}(X)\to \mathbb{R} $ defined by
	$$h(A, B) = \max \{d(A,B), d(B,A)\}$$ is a metric on $\mathbbm{H}(X)$ called the\textit{ Hausdorff-Pompeiu metric} or \textit{ Hausdorff metric}.
\end{definition}
The following is a property of the Hausdorff metric.
\begin{theorem}\cite{barnsley2014fractals}
	Let $(X,d)$ be a metric space. Then $$h(A\cup B, C\cup D)\le\max\{h(A, C), h(B,D)\}$$ for every $A,B,C,D \in \mathbbm{H}(X). $
\end{theorem}
%
An iterated function system is a dynamical system with finitely many contraction maps acting on a complete metric space. The following is the formal definition of a hyperbolic $IFS$.
\begin{definition}\cite{barnsley2014fractals}
	A \textit{hyperbolic Iterated Function System }(hyperbolic IFS) consists of a complete metric space $(X,d)$ and a finite number of contraction mappings $f_i:X\to X$, with respective contractivity factors $t_i$ for $i=1,2,\cdots,n.$
	
\end{definition}
The $IFS$ is a method to construct objects with fractal nature mathematically. The attractor or the set fixed point obtained from the $IFSs$ are mostly self-similar fractals. We have the following theorem for the existence and uniqueness of the attractor of an $IFS$.
\begin{theorem}\cite{barnsley2014fractals}
	Let $\displaystyle \Bigl\{X;f_i,i=1,2,\cdots,n\Bigr\}$ be an IFS with contractivity factor $t$. Then the transformation $W:\mathbbm{H}(X) \to \mathbbm{H}(X)$ defined by $\displaystyle W(B)=\cup_{i=1}^nf_i(B)$ for all $B\in \mathbbm{H}(X)$, is a contraction mapping on the complete metric space $(\mathbbm{H}(X),h)$ with contractivity factor $t$. Its unique fixed point, $A\in \mathbbm{H}(X)$, exists and is given by $\displaystyle A=\lim_{n\to \infty}W^{[n]}(B)$ for any $B\in \mathbbm{H}(X)$.
\end{theorem}
To construct a given object with fractal nature mathematically one can use the aid of $IFS$ theory. The following theorem called the \textit{Collage theorem} ensures how this can be achieved.
\begin{theorem}\cite{barnsley2014fractals}
	Let $L\in \mathbbm{H}(X)$ and $\epsilon>0$ be given. Choose an $IFS$ $\displaystyle\Big\{X;f_i,i=1,2,\cdots,n\Bigr\}$ with contractivity factor $0\le t<1$, so that $h\Big(L,\displaystyle\bigcup_{i=1}^{n}f_i(L)\Big)\le \epsilon$
	Then $h\Big(L,F\displaystyle\Big)\le \frac{\epsilon}{1-t},$
	where $F$ is the attractor of the IFS.
\end{theorem}

There have been considerable developments in the $IFS$ theory, and many generalizations to this have come in the recent years. A generalization in which the metric space is replaced with a partial metric space was done in \cite{minirani2014fractals}.
The partial metric is defined as follows.
\begin{definition}\cite{matthews1994partial}
	A partial metric space is a pair $(X,p:X\times X\to \mathbb{R})$ such that
	\begin{itemize}
		\item[P1:] $0\le p(x,x)\le p(x,y)$ (non-negativity and small self distances),
		\item[P2:] if $p(x,x)=p(x,y)=p(y,y)$ then $x=y$ (indistancy implies equality),
		\item[P3:] $p(x,y)=p(y,x)$ (symmetry), and
		\item[P4:] $p(x,z)\le p(x,y)+p(y,z)-p(y,y)$ (triangularity)
	\end{itemize}
\end{definition}
The class of functions such as Lipschitz functions and contraction functions in partial metric space can be defined in a similar manner to the metric space theory. The contraction in a partial metric space is defined as follows.
\begin{definition}\cite{matthews1994partial}
	For each partial metric space $(X,p)$, a contraction is a function $f:X\to X$ for which there exists a $c\in [0,1)$ such that for all $x,y\in X$, $p(f(x),f(y))\le c\cdot p(x,y)$.
\end{definition}
The Banach contraction principle, which ensures the existence and uniqueness of fixed point for a contraction mapping in a complete metric space, is central to the fixed point theory and it makes the $IFS$ theory possible. In \cite{matthews1994partial}, Matthews generalized the Banach contraction principle to the partial metric spaces.
The statement of the theorem is as follows.
\begin{theorem}\cite{matthews1994partial}
	For each contraction $f$ over a complete partial metric space $(X,p)$ there exists a unique $x\in X$ such that $x=f(x)$. Also $p(x,x)=0.$
\end{theorem}
In \cite{minirani2014fractals}, Minirani considered $IFS$ with partial metric space in the system and termed it as partial iterated function system($IFS_p$).
The following is the definition of a hyperbolic $IFS_p$.
\begin{definition}\cite{minirani2014fractals}
	A (hyperbolic) partial iterated function system($IFS_p$) consists of a complete partial metric space $(X,p)$ together with a finite set of contraction mappings $w_j:X\to X$ with respective contractivity factors $s_j$, for $j=1,2,\cdots,N.$ Further, $s=\max\{s_j:j=1,2,\cdots, N \}$ is called the contraction factor of the $IFS_p$.
\end{definition}
A homogeneous $IFS_p$ is defined as follows.
\begin{definition}\cite{minirani2014fractals}
An $IFS_p$ is said to be homogeneous if the contraction factor $s_j=s$ for all $j=1,2,\cdots,N.$
\end{definition}
The following theorem provides the existence and uniqueness of the attractor of an $IFS_p$.
\begin{theorem}\cite{minirani2014fractals}
Let $\{X;w_j, j=1,2,\cdots,N\}$ be a hyperbolic $IFS_p$ with contraction factor $s$. Then the transformation $W: \mathbbm{H}_p(X)\to \mathbbm{H}_p(X)$ defined by $W(B)=\cup_{i=1}^N w_j(B)$ for all $B\in \mathbbm{H}_p(X)$ is a contraction mapping on the complete metric space $(\mathbbm{H}_p(X),h_p)$ with contraction factor $s$. The unique fixed point $A\in \mathbbm{H}_p(X)$, called the attractor of this $IFS_p$, obeys $A=W(A)=\cup_{i=1}^Nw_j(A)$ and is given by $A=\lim_{n\to \infty} W^{\circ n}(B)$ for any $B\in \mathbbm{H}_p(X)$.
\end{theorem}
Recently, study in $IFS$ theory was connected to algebra and several results have been obtained. In the book \cite{barnsley2006superfractals}, Barnsley defines the $IFS$ semigroup as follows.
\begin{definition}\cite{barnsley2006superfractals}
	An $IFS$ semigroup is a semigroup of transformations generated by an $IFS$. For an $IFS$, $\{X;f_1,f_2,\cdots,f_N\}$ the corresponding $IFS$ semigroup is denoted by $\mathcal{S}_{\{X;f_1,f_2,\cdots,f_N \}}$.
\end{definition}	
In the next section we provide our definitions and major results. We prove the continuity of a partial metric function. Also, we prove the collage theorem for a partial iterated function system. A method to address the attractor of an $IFS_p$ is also obtained. The continuity of the map which associates to each contraction in a complete partial metric space, the corresponding fixed point is obtained. Also, $IFS_p$ semigroup is defined and related results are obtained.
\section{Definitions and Major Results}
In this section we discuss certain results in partial metric spaces. We are proving the continuity of the partial metric function. Also, results leading to Collage theorem in $IFS_p$ is provided. A method to address the points in the attractor of an $IFS_p$ is given, and certain results are obtained using this addressing method. Later we discuss certain results concerning the space of contractions on a complete partial metric space with fixed contractivity factor. Also, continuity of the map is obtained, that associates to each contraction in a complete partial metric space the corresponding unique fixed point. Afterwards, we are defining the $IFS_p$ semigroup and proving certain related results. 

The following theorem provides a proof for the continuity of the partial metric function. 
\begin{theorem}
	Let $(X,p)$ be a partial metric space. Let $a\in X$, then $p(a,\cdot):X\to X$ is continuous.
	\begin{proof}
		We have	from the properties of partial metric \begin{align*}
		p(a,x_1)&\le p(a,x_2)+p(x_1,x_2)-p(x_2,x_2) \quad\text{ and} \\
		p(a,x_2)&\le p(a,x_1)+p(x_1,x_2)-p(x_1,x_1) \\
		\implies	p(a,x_1)-p(a,x_2)&\le p(x_1,x_2)-p(x_2,x_2) \quad\text{ and} \\
		p(a,x_2)-p(a,x_1)&\le p(x_1,x_2)-p(x_1,x_1) \\
		\implies |p(a,x_1)-p(a,x_2)|&\le p(x_1,x_2)-\min\{p(x_1,x_1),p(x_2,x_2)\}\le p(x_1,x_2) 
		\end{align*}
		Hence the partial metric function $p(a,\cdot)$  is continuous. 
	\end{proof}
\end{theorem}
Now we present a theorem which will help us in proving the collage theorem for the partial iterated function systems.
\begin{theorem}\label{col}
	Let $(X,p)$ be a complete partial metric space, $f:X\to X$ be a contraction mapping with contractivity factor $0\le s<1$ and let the fixed point of $f$ be $x_f\in X$. Then $p(x,x_f)\le \frac{1}{1-s}p(x,f(x))$. 
	\begin{proof}
		The partial metric function $p(a,x)$ for fixed $a\in X$ is continuous in $x\in X$. Using the properties of the partial metric, we have
		\begin{align*}
		p(x,x_f)&=p(x,\lim_{n\to \infty}f^n(x))\\
		&=\lim_{n\to \infty} p(x,f^n(x))
		\qquad[\text{by continuity of } p(x,.)] \\
		&\le \lim_{n\to \infty}\Big[\sum_{m=0}^n p(f^{m-1}(x),f^m(x))-\sum_{m=1}^{n-1} p(f^{m}(x),f^m(x))\Big]\\
		&\le \lim_{n\to \infty}\sum_{m=0}^n p(f^{m-1}(x),f^m(x))\\
		&\le \lim_{n\to \infty} p(x,f(x))\sum_{m=0}^n s^m\\
		&\le \frac{1}{1-s} p(x,f(x)).   	
		\end{align*}
	\end{proof}
\end{theorem}
Now, using the usual definition of condensation transformation we define the family $\mathcal{C}$ as follows:
\begin{definition}
	Let $(X,p)$ be a partial metric space and let $C\in \mathbbm{H}_p(X)$. Define $w_0:\mathbbm{H}_p(X)\to \mathbbm{H}_p(X)$ by $w_0(B)=C, \quad \forall B \in \mathbbm{H}_p(X)$. Then $w_0$ is called a condensation transformation and $C$ is called the associated condensation set.\\
	Let $\mathcal{C}:=\{C\in \mathbbm{H}_p(X)/h_p(C,C)=0\}$. 
\end{definition}
Even though the condensation transformation in $IFS$ theory are contractions, this is not the case in $IFS_p$. The following example shows that a condensation transformation need not be a contraction in the partial metric case.
\begin{example}
	Consider $X=\mathbb{R}^+\cup\{0\}$ with the partial metric $p(x,y)=\max(x,y)$. Let  $C=[3,4]$, then $C\in \mathbbm{H}_p(X)$. Now let $w_0: \mathbbm{H}_p(X)\to  \mathbbm{H}_p(X)$ be the condensation $w_0(B)=C, \quad\forall B\in  \mathbbm{H}_p(X)$.\\
	Let $B_1=[0,1]$ and $B_2=[2,3]$.\\
	Then \begin{align*}h_p(w_0(B_1),w_0(B_2))&=h_p(C,C)\\&=\sup\{p(x,x):x\in C\}\\&=\sup\{\max{x,x}:x\in C\}\\&=\sup\{x:x\in C \}=4.\end{align*} 
	Now, \begin{align*}h_p(B_1,B_2)&=\max\{\rho_p(B_1,B_2),\rho_p(B_2,B_1) \}\\&=\max\{\sup_{x\in[0,1]}p(x,[2,3]),\sup_{x\in[2,3]}p(x,[0,1])\}\\&=\max\{\sup_{x\in[0,1]}2,\sup_{x\in[2,3]}x\}\\&=\max\{2,3\}=3.\end{align*} Hence $h_p(w_0(B_1),w_0(B_2))=4>3=h_p(B_1,B_2).$ Thus, the condensation transformation $w_0$ is not a contraction in the partial metric space $(X,p)$.
\end{example}

In the light of the above example we obtain the following theorem.
\begin{theorem}
	The condensation transformation is a contraction if and only if the corresponding condensation set $C\in\mathcal{C}$.
	\begin{proof}
		Let $w_0$ be a condensation transformation with the associated condensation set $C$.\\	Suppose $C\in \mathcal{C}$, then $h_p(w_0(B_1), w_0(B_2))=h_p(C,C)=0, \forall B_1,B_2\in \mathbbm{H}_p(X)$. Since, $h_p(B_1, B_2)\ge 0$, we get $h_p(w_0(B_1), w_0(B_2))\le h_p(B_1, B_2)$. Hence $w_0$ is a contraction.\\Now we prove the   contrapositive of the converse part.\\For, suppose $h_p(C,C)\ne 0$. To prove $w_0$ is not a contraction. Choose $B_1=B_2=C$, then $h_p(w_0(B_1),w_0(B_2))=h_p(C,C)\nless h_p(C,C)=h_p(B_1,B_2).$ Hence $w_0$ is not a contraction.
	\end{proof}
\end{theorem}
The above theorem explains why we need the following additional condition while defining the Condensation $IFS_p$.
Now we define the condensation $IFS_p$.
\begin{definition}[Condensation $IFS_p$] An $IFS_p$ is a condensation $IFS_p$ if the condensation set $C$ corresponding to the condensation transformation satisfies $h_p(C,C)=0$, $i.e.,$ if $C\in \mathcal{C}$.
	
\end{definition}
Now we state the Collage theorem in partial metric space, the proof of which follows from the Theorem \ref{col}
.\begin{theorem}[Collage Theorem in Partial Metric Space]
	Let $(X,p)$ be a complete partial metric space. Let $L\in \mathbbm{H}_p(X)$  and let $\epsilon\ge0$ be given. Choose an $IFS_p$ ( or $IFS_p$ with condensation) $\{X;(w_0),w_1,w_2,\cdots,w_n\}$ with contractivity factor $0\le s< 1$ so that, $h_p(L,\cup_{m=1\\(m=0)}^n w_m(L))\le \epsilon$ where $h_p$ is the Hausdorff partial metric. Then $h_p(L,A)\le \frac{\epsilon}{1-s}$ where $A$ is the attractor of the $IFS_p$. Equivalently, $h_p(L,A)\le (1-s)^{-1}h_p(L,\cup_{m=1\\(m=0)}^n w_m(L)), \forall L\in \mathbbm{H}_p(X)$.
\end{theorem}

Further we introduce the \textit{shift space} for a partial IFS.

Let $I=\{1,2,3,\cdots,n\}$ and $\Lambda=\Lambda(I)$ be the space of infinite words with letters from $I$. We define $p_{\Lambda}:\Lambda\times \Lambda\to \mathbb{R}$ by, 
$$p_{\Lambda}(\alpha,\beta)=\sum_{k=1}^\infty\frac{1-\delta_{\alpha_k}^{\beta_k}}{3^k}+\frac{1}{k^2},$$ where $\delta_{\alpha_k}^{\beta_k}$ is the Kronecker-$\delta$ function.
Then $(\Lambda,p_{\Lambda})$ is a partial metric space. 
\begin{definition}
	The partial metric space $(\Lambda(I),p_{\Lambda})$ is called the shift space associated with the partial IFS $S=(X,(f_i)_{i\in I}).$
\end{definition}
Now we provide a result concerning the location of the unique fixed point for a contraction in a complete partial metric space provided the contraction satisfies certain condition. 
\begin{theorem}\label{fix}
	Let $(X,p)$ be a complete partial metric space. Let $f:X\to X$ be a contraction and $H$ be a closed set such that $f(H)\subseteq H.$ If $x_f$ is the fixed point of $f$, then $x_f\in H.$
	\begin{proof}
		Let $a_0\in H$ and let $(a_n)_{n\in \mathbb{N}}$ be the sequence defined by $a_{n+1}=f(a_n),\forall n\in \mathbb{N}$. Then since $a_n\in H,\forall n\in \mathbb{N},$ $\displaystyle \lim_{n\to \infty}a_n=x_f$ and $H$ is closed, we get $x_f\in H$.
	\end{proof}
\end{theorem}
We generalize the idea of a Lipschitz function and a contraction to partial metric space as follows.
\begin{definition}
	Let $(X,p)$ be a partial metric space. For a transformation $f:X\to X$ we define the Lipschitz constant associated with $f$ as $$\lipp(f)= \sup_{\substack{x,y\in X\\x\ne y}}\frac{p(f(x),f(y))}{p(x,y)}\in[0,+\infty].$$ If $\lipp(f)<+\infty$ then $f$ is said to be a Lipschitz function, and if $\lipp(f)<1$ then $f$ is called a contraction.
\end{definition}
Now we prove a result related to Lipschitz constant of a function in a partial metric space. This result discuss how the diameter of a set varies with the Lipschitz constant of a transformation applied to the set. 
\begin{theorem}\label{diam}
	For a partial metric space $(X,p)$ and a function $f:X\to X$ the following result holds for every subset $A$ of $X$:$$\diam(f(A))\le \lipp(f)\diam(A).$$
	\begin{proof}
		We have \begin{align*} \diam(f(A))&=\sup_{x,y\in A}p(f(x),f(y))\\
		&\le sup_{x,y\in A}\lipp(f)p(x,y)\\
		&=\lipp(f)\sup_{x,y\in A} p(x,y)\\
		\therefore \diam(f(A))&\le \lipp(f)\diam(A).
		\end{align*}
	\end{proof}
\end{theorem} 
The following results are related to addressing the attractor of an $IFS_p$.
\begin{theorem}
	Let $S=(X,(f_i)_{i\in I})$ be a partial IFS, where $(X,p)$ is a complete partial metric space. Also, let $A=A(S)$ be the attractor of $S$ and $c=\sup_{i\in I}\lipp(f_i)<1$. Then the following are true:\begin{itemize}
		\item[(i)] For $m\in \mathbb{N},$ we have $A_{[w]_{m+1}}\subseteq A_{[w]_{m}}$, for all $w\in\Lambda=\Lambda(I)$, and $\lim_{m\to \infty}\diam(A_{[w]_m})=0$, where $A_{[w]_m}=f_{[w]_m}(A)=f_{w_1}\circ f_{w_2}\circ \cdots f_{w_m}(A)$.
		\item[(ii)] If $a_w$ is defined by $\{a_w\}=\cap_{m\in \mathbb{N}^*}\overline{A_{[w]_m}}$, then $\lim_{m\to \infty}p(x_{[w]_m},a_w)=0$, where $x_{[w]_m}$ denotes the fixed point of the contraction $f=f_{[w]_m}=f_{w_1}\circ f_{w_2}\circ \cdots f_{w_m}$.
		\item[(iii)] Foe every $a\in A$ and every $w\in \Lambda$ we have $\lim_{m\to \infty}f_{[w]_m}(a)=a_w.$
	\end{itemize}
	\begin{proof}
		\begin{enumerate}[label=(\roman*)]
			\item Let $w\in \Lambda$ be given. We first prove that $A_{[w]_{m+1}}\subseteq A_{[w]_{m}}$, for all $m\in \mathbb{N}.$ Since $A$ being the attractor of $S$, we have $A=\cup_{i=1}^nf_i(A)$, and hence we get $A_{[w]_1}=f_{w_1}(A)\subseteq A.$ Now for every $m\in \mathbb{N},m\ge2$, we have, $A_{[w]_m}=f_{[w]_m}(A)=f_{[w]_{m-1}}\circ f_{w_m}(A)\subseteq f_{[w]_{m-1}}(A)=A_{[w]_{m-1}}$. From the \textbf{\small Theorem \ref{diam}.} we have $\diam(A_{[w]_m})\le c^m\diam(A),\forall m\in \mathbb{N}.$
		
			\item From (i) since $\overline{A_{[w]_{m+1}}}\subseteq \overline{A_{[w]_m}}$ and since $\diam({A_{[w]_m}})\to 0$ we get $\cap_{m\in \mathbb{N}^*}\overline{A_{[w]_m}}$ is a singleton set.
			Let $\cap_{m\in \mathbb{N}^*}\overline{A_{[w]_m}}=\{a_w\}.$ Since $f_{[w]_m}$ being a contraction is continuous and $\overline{A_{[w]_m}}\subseteq A$, we have for every $m\in \mathbb{N}^*$, $f_{[w]_m}(\overline{A_{[w]_m}})\subseteq \overline{f_{[w]_m}(A)}=\overline{A_{[w]_m}}.$ Now by \textbf{\small Theorem \ref{fix}} , we have $x_{[w]_m}\in \overline{A_{[w]_m}}$. Then 
			\begin{align*}
			p(x_{[w]_m},a_w)&\le \diam(\overline{A_{[w]_m}})\\
			&\le \diam({A_{[w]_m}})\\
			&\le c^m \diam(A),
			\end{align*}
			for every $m\in \mathbb{N}^*$ and as a consequence we have $\displaystyle\lim_{m\to \infty}p(x_{[w]_m},a_w)=0$.
			\item[(iii)] We have 	\begin{align*}
			p(f_{[w]_m}(a),a_w)&\le \diam(\overline{A_{[w]_m}})\\
			&\le \diam({A_{[w]_m}})\\
			&\le c^m \diam(A),
			\end{align*}
			for every $m\in \mathbb{N}^*$ and $\lim_{m\to \infty}c^m=0.$ Hence for every $a\in A$ and $w\in \Lambda$, we have  $\displaystyle\lim_{m\to \infty}f_{[w]_m}(a)=a_w.$
		\end{enumerate}
	\end{proof}
\end{theorem}
 Now we introduce a partial metric on the set of contractions in a partial metric space.
\begin{definition}
	Let $(X,p)$ be a complete partial metric space. We denote the set of all contractions with contractivity factor $t$ on $X$ as $Con_t(X,p)$, i.e., $Con_t(X,p)=\{f:X\to X| p(f(x),f(y))\le t\cdot p(x,y) \}$, where $t\in [0,1)$. We define a partial metric on $Con_t(X,p)$ as $\displaystyle\bar{p}(f,g):=\min\{1, \sup_{x\in X} p(f(x),g(x)) \}$.
\end{definition}
The following theorem proves the completeness of the space $(Con_t(X,p),\bar{p})$ provided $(X,p)$ is complete. 
\begin{theorem}
	Let $(X,p)$ be a complete partial metric space. Then $(Con_t(X,p),\bar{p})$ is a complete partial metric space.
	\begin{proof}
		Let $(f_n)_{n=1}^{\infty}$ be a Cauchy sequence in $(Con_t(X,p),\bar{p})$. Then we have $\bar{p}(f_m,f_n)\to k$ for some $k\in [0,\infty)$ as $m,n\to \infty$, i.e.,\\
		$\displaystyle\lim_{m,n\to \infty} \bar{p}(f_m,f_n)=\lim_{m,n\to \infty}\min\{1,\sup_{x\in X}p(f_m(x),f_n(x)) \}=k.$ \\We observe that $k\le 1$.\\
		We first prove that $(f_n(x))$ is a Cauchy sequence in $(X,p)$, for each $x\in X$. For, it is enough to show that $\displaystyle \lim_{m,n\to \infty}p(f_m(x),f_n(x))=l$, for some $l\in [0,\infty)$.
		Suppose for some $x\in X$, $\displaystyle \lim_{m,n\to\infty}p(f_m(x),f_n(x))$ does not exist. Then by definition of $\bar{p}$ we get, $\displaystyle \lim_{m,n\to\infty}\bar{p}(f_m,f_n)$ does not exist, which is a contradiction to the assumption that $(f_n)_{n=1}^{\infty}$ is a Cauchy sequence.\\
		Hence $(f_n(x))$ is a Cauchy sequence in $(X,p)$, for each $x\in X$.\\
		We define $f_n(x):=\lim_{m,n\to \infty}p(f_m(x),f_n(x) )$ for each $x\in X$.\\
		Now we will prove that $f\in Con_t(X,p)$.\\
		For, let $x,y\in X$. Then \begin{align*}
		p(f(x),f(y))&=p(\lim_{n\to \infty}f_n(x),\lim_{n\to \infty}f_n(y))\\
		& =\lim_{n\to \infty}p(f_n(x),f_n(y)) \quad [\because \text{ partial metric is a continuous function}]\\
		& \le \lim_{n\to \infty}t\cdot p(x,y)\quad[\because f_n'\text{s are Cauchy sequences} ]\\
		\text{i.e., } p(f(x),f(y))&\le t\cdot p(x,y), \, \forall x,y\in X.
		\end{align*}
		Thus $f\in Con_t(X,p)$.\\ Hence every Cauchy sequence converges in the partial metric space $(Con_t(X,p),\bar{p})$, which proves that $(Con_t(X,p),\bar{p})$ is a complete partial metric space.
	\end{proof}
\end{theorem}
In the next theorem we prove the continuity of the map which associates to each contraction in a complete partial metric space, the corresponding fixed point.
\begin{theorem}
	Let $(X,p)$ be a complete partial metric space and let $r: Con_t(X,p)\to (X,p)$ be defined by $r(f)=x_f$ for each $f\in Con_t(X)$ where $x_f$ denotes the fixed point of the contraction $f$. Then $r$ is continuous.
	\begin{proof}
		To prove $r$ is continuous, it is enough to show that  $r(f_n)\to r(f)$ in $(X,p)$ for every $(f_n)\to f$ in $(Con_t(X,p),\bar{p})$.\\
		For we have,
		\begin{align*}
		\lim_{n\to \infty}p(r(f_n),r(f))&=\lim_{n\to \infty}p(x_{f_n},x_f)\\
		&=\lim_{n\to \infty}p(\lim_{m\to \infty}f_n^m(x),\lim_{m\to \infty}f^m(x)), \quad \text{ for some } x\in X\\
		&=\lim_{n\to \infty}\lim_{m\to \infty}p(f_n^m(x),f^m(x)) \quad [\text{by continuity of } p] \\
		&= \lim_{m\to \infty}p(\lim_{n\to \infty}f_n^m(x),f^m(x)) \quad [\text{by continuity of } p] \\
		&= \lim_{m\to \infty}p(f^m(x),f^m(x)) \quad [f_n \text{ converges to  } f] \\
		&=p(\lim_{m\to \infty}f^m(x),\lim_{m\to \infty}f^m(x)) \quad [\text{ by continuity of } p]\\
		&=p(x_f,x_f)\\
		&=p(r(f),r(f))
		\end{align*}
		Thus, $r(f_n)\to r(f)$ in $(X,p)$ for every $(f_n)\to f$ in $(Con_t(X,p),\bar{p})$.\\
		Hence $r: Con_t(X,p)\to (X,p)$ is continuous.
	\end{proof}
\end{theorem}
Now we provide some algebraic results in $IFS_p$. We define the $IFS_p$ semigroup as follows.
\begin{definition}
	An $IFS_p$ semigroup is a semigroup of transformations generated by an $IFS_p$. For an $IFS_p$, $\{(X,p); f_1,f_2,\cdots,f_N \}$, the $IFS_p$ semigroup will be denoted by $\mathcal{S}_{\{(X,p);f_1,f_2,\cdots,f_N\}}^p$.
\end{definition}
The following theorem proves that the set of Lipschitz transformations on a partial metric space with function composition operation is a semigroup.
\begin{theorem}
	Let $(X,p)$ be a partial metric space. Then, the set of Lipschitz transformations on $X$ forms a semigroup.
	\begin{proof}
		Let $\mathcal{S}=\{f:X\to X \,|\,\lipp(f)<\infty \}$.\\
		To show $(\mathcal{S},\circ)$ is closed and associative.\\
		Let $f_1,f_2\in \mathcal{S}$. Then $\lipp(f_1)<\infty$ and $\lipp(f_2)<\infty$.\\
		We have \begin{align*}
		\lipp(f_1\circ f_2) &=\sup_{x,y\in X;x\ne y}\frac{p(f_1\circ f_2(x),f_1\circ f_2(y))}{p(x,y)}\\
		&=\sup_{x,y\in X;x\ne y}\frac{p(f_1 (f_2(x)),f_1 (f_2(y)))}{p(x,y)}\\
		&\le \sup_{x,y\in X;x\ne y}\frac{p(f_1(x),f_1(y))}{p(x,y)}\\
		&=\lipp(f_1)<\infty
		\end{align*}
		Thus $f_1\circ f_2\in \mathcal{S},\, \forall f_1,f_2\in \mathcal{S}$ and hence $(\mathcal{S},\circ)$ is closed.\\
		Also we have $$(f_1\circ f_2)\circ f_3(x)=f_1 (f_2( f_3(x))=f_1\circ(f_2\circ f_3)(x), \forall x\in X$$
		Thus, $(f_1\circ f_2)\circ f_3=f_1\circ(f_2\circ f_3)$.\\
		Hence $(\mathcal{S},\circ)$ is associative.\\
		Therefore, $(\mathcal{S},\circ)$ is a semigroup.
		
	\end{proof}
\end{theorem}
The next theorem proves that the set of contractions on a partial metric space with function composition operation is a semigroup.
\begin{theorem}
	Let $(X,p)$ be a partial metric space. Then, the set of contractions on $X$ 
	forms a semigroup.
	\begin{proof}
		We have the contractions on $X$ are the Lipschitz transformations on $X$ with Lipschitz constant less than 1.\\
		Let $\mathcal{S}=\{f:X\to X \,|\,\lipp(f)<1 \}$.\\
		To show $(\mathcal{S},\circ)$ is closed and associative.\\
		Let $f_1,f_2\in \mathcal{S}$. Then $\lipp(f_1)<1$ and $\lipp(f_2)<1$.\\
		We have \begin{align*}
		\lipp(f_1\circ f_2) &=\sup_{x,y\in X;x\ne y}\frac{p(f_1\circ f_2(x),f_1\circ f_2(y))}{p(x,y)}\\
		&=\sup_{x,y\in X;x\ne y}\frac{p(f_1 (f_2(x)),f_1 (f_2(y)))}{p(x,y)}\\
		&\le \sup_{x,y\in X;x\ne y}\frac{p(f_1(x),f_1(y))}{p(x,y)}\\
		&=\lipp(f_1)<1
		\end{align*}
		Thus $f_1\circ f_2\in \mathcal{S}$ and hence $(\mathcal{S},\circ)$ is closed.\\
		The associative property of $\circ$ in $\mathcal{S}$ follows as a hereditary property.  \\
		Therefore, $(\mathcal{S},\circ)$ is a semigroup.
	\end{proof}
\end{theorem}
We now provide an example to show that the set of all Lipschitz functions with a fixed Lipschitz constant $L>1$ does not in general form a semigroup.
\begin{example}
	Consider $X=[0,1]$ with the partial metric $p(x,y)=\max\{x,y \}, \forall x,y\in X$. Let $\mathcal{S}=\{f:[0,1]\to[0,1]| \lipp(f)=2 \}$ and let $f(x)=2x$ and $g(x)=2x^2$. Then, we have
	\begin{align*}
	&\lipp(f)=\sup_{\substack{x,y\in [0,1]\\x\ne y}}\frac{p(f(x),f(y))}{p(x,y)}
	=\sup_{\substack{x,y\in [0,1]\\x\ne y}}\frac{p(2x,2y)}{p(x,y)}
	=\sup_{\substack{x,y\in [0,1]\\x\ne y}}\frac{\max\{2x,2y\}}{\max\{x,y\}}
	=2, \quad \text{ and }
	\\
	&\lipp(g)=\sup_{\substack{x,y\in [0,1]\\x\ne y}}\frac{p(g(x),g(y))}{p(x,y)}
	=\sup_{\substack{x,y\in [0,1]\\x\ne y}}\frac{p(2x^2,2y^2)}{p(x,y)}
	=\sup_{\substack{x,y\in [0,1]\\x\ne y}}\frac{\max\{2x^2,2y^2\}}{\max\{x,y\}}
	=2.
	\end{align*}
	Therefore, $f,g\in \mathcal{S}$. But, we have
	\begin{equation*}
	\lipp(f\circ g)
	=\sup_{\substack{x,y\in [0,1]\\x\ne y}}\frac{p(f\circ g(x),f\circ g(y))}{p(x,y)}
	=\sup_{\substack{x,y\in [0,1]\\x\ne y}}\frac{p(4x^2,4y^2)}{p(x,y)}
	=\sup_{\substack{x,y\in [0,1]\\x\ne y}}\frac{\max\{4x^2,4y^2\}}{\max\{x,y\}}
	=4.
	\end{equation*}
	Hence $f\circ g \notin \mathcal{S}.$\\
	Therefore, $(\mathcal{S},\circ)$ is not closed and hence not a semigroup.
\end{example}

We conclude the findings of our paper in the next section.
\section{Conclusion}
The idea of partial metric spaces was developed to incorporate the metric notion in those spaces that are not Hausdorff. The concept of fractals to such non-Hausdorff spaces can be established using the definition of partial iterated function systems. This article proved certain results on partial iterated function system($IFS_p$). Also, it generalized some concepts from classical $IFS$ theory to the theory of $IFS_p$. Further, it proved the Collage theorem for partial iterated function system.

Introducing a proper shift space to address the points in the attractor of an $IFS$ is critical to understand the dynamics and the separation properties of it. This paper introduced a proper shift space for the $IFS_p$ to address its attractor. It also discusses results related to the address of points in the attractor of an $IFS_p$. 

This article considered the partial metric space of contractions with a fixed contractivity factor and proved the completeness of it, under suitable conditions. The continuity of the map which associates to each contraction in a complete partial metric space, its corresponding fixed point has been obtained. The paper also defined the $IFS_p$ semigroup and proved that under function composition, the set of Lipschitz maps and the set of contractions are semigroups.
\section*{Acknowledgement}
The first author is very grateful to Council of Scientific \& Industrial Research(CSIR), India for their financial support.

\bibliographystyle{unsrt}
\bibliography{reference}{}
\end{document}